\documentclass{amsart}
\usepackage{amssymb}

\newtheorem{theorem}{Theorem}[section]
\newtheorem{lemma}[theorem]{Lemma}

\newtheorem{proposition}[theorem]{Proposition}

\theoremstyle{definition}

\theoremstyle{remark}

\begin{document}

\title[On stability of classes of solutions]
{On stability of classes of solutions to partial differential relations 
constructed by quasiconvex functions and null Lagrangians 
with respect to precompact famalies in $C_{loc}$}

\author[A.~A.~Egorov]{A.~A.~Egorov$^1$}
\address{Sobolev Institute of Mathematics, Novosibirsk, Russia}
\email{yegorov@math.nsc.ru}

\begin{abstract} 
We prove theorems on stability of classes of solutions to partial differential relations 
constructed by quasiconvex functions and null Lagrangians
with respect to precompact famalies in $C_{\operatorname{loc}}$.
\end{abstract}

\maketitle

\addtocounter{footnote}{+1}\footnotetext{The author is supported by 
the Russian Foundation for Basic Research grant No. 20--01--00661.}

Let
$\mathfrak G$ 
be the class of
$W_{\operatorname{loc}}^{1,k}$-solutions 
$u\colon V\to\mathbb R^m$
(defined on domains 
$V\subset\mathbb R^n$)
to the equation
\begin{equation}
\label{eq:osnovnoe-uravnenie}
F(u'(x))=G(u'(x))\quad\text{a.e. }x\in V,
\end{equation}
where
$F\colon \mathbb R^{m\times n}\to\mathbb R$ 
is a nonnegative quasiconvex function
and
$G\colon \mathbb R^{m\times n}\to\mathbb R$ 
is a null Lagrangian.
Here 
$u'(x)$
denotes the Jacobi matrix 
of~$u$
at 
$x\in V$.

Let 
$\mathfrak F$
be the class of mappings 
$v\in W_{\operatorname{loc}}^{1,k}(V;\mathbb R^m)$ 
(defined on domains
$V\subset\mathbb R^n$) 
for which there exists a finite measurable function 
$K:V\to[1,+\infty)$, 
finite almost everywhere, such that
$$
F(v'(x))\le K(x)G(v'(x))\quad\text{a.e. }V.
$$
Then for 
$v\colon V\to\mathbb R^m$
of the 
class~$\mathfrak F$
and for a.e.
$x\in V$
we can define 
$$
K(x,v)
=\begin{cases} 
\frac{F(v'(x))}{G(v'(x))}&\text{if }G(v'(x))>0;\\
1&\text{if }F(v'(x))=0.
\end{cases}
$$

The 
class~$\mathfrak G$
has some stability property if any mapping 
$v\in\mathfrak F$ 
for which the function
$K(x,v)$
is close 
to~$1$
also close to some mapping
$u\in\mathfrak G$.

In~\cite{Egor2008}, 
the author has obtained some results on stability 
of~$\mathfrak G$
in the case when the discrepancy between 
$K(x,v)$
and~$1$
is measured in the norm of
$L^\infty(V)$.
In this case 
$v$
belongs to the classes
$\mathfrak G(K):=\{v\colon V\to\mathbb R^m,\ v\in\mathfrak F,\ \operatorname{ess\,sup}_{x\in V}K(x,V)<K\}$,
$K\ge1$.
Note that the class 
$\mathfrak G(K)$
consists of
$W_{\operatorname{loc}}^{1,k}$-solutions 
$v\colon V\to\mathbb R^m$
(defined on domains 
$V\subset\mathbb R^n$)
of the inequality
\begin{equation}
\label{eq:nl-neravenstvo}
F(v'(x))\le KG(v'(x))\quad\text{a.e.\ }V.
\end{equation}

The aim of the present paper is to prove that a mapping 
$v\in\mathfrak G(K)$
is close to some 
$u\in\mathfrak G$
in the case when the function 
$K(x,v)$
is close 
to~$1$
only some integral sense.

Our results are analogues of N.~A.~Kudryavtseva and Yu.~G.~Reshetnyak's 
results~\cite{KudrR1993} 
on stability of M\"obius transformations 
with respect to precompact (in $C_{\operatorname{loc}}$) famalies 
of mappings with bounded distortion.
A mapping 
$v\in W_{\operatorname{loc}}^{1,n}(V;\Bbb R^n)$
of an open set
$V\subset\Bbb R^n$ 
is an (orientation-preserving) mapping with 
$K$-bounded 
distortion,
$K\ge1$,
if 
$v$
satisfies the distortion inequality
\begin{equation}
\label{eq:diltn}
|v'(x)|^n\le K\det v'(x)\quad\text{a.e.\ }V,
\end{equation}
where
$|v'(x)|$
is the operator norm of the matrix
$v'(x)$.
If, in addition,
$v$
is topological, then 
$v$
is 
$K$-quasiconformal.
The distortion inequality is the particular case 
of~\eqref{eq:nl-neravenstvo}
with the following functions
$F(v'(x))=|v'(x)|^n$
and 
$G(v'(x))=\det v'(x)$.
The theory of quasiconformal mappings and mappings with bounded 
distortion is the key part of modern geometric analysis which 
has many diverse applications, for example, see monographs 
\cite{AstaIM2009,BojaGMR2013,Beli1974,GoldR1983,GoldR1990,
GutlRSY2012,IwanM2001,Kopy1990,Krus1975,Krus1979,LehtV1973,
Mart2000,MartRSY2009,Resh1982a,Resh1982b,Resh1989,Resh1994,
Resh1996,Rick1993,Vais1971,Vour1988}
and the bibliography therein.
In this monographs the results on stability of M\"obius transformations
are playing an important role.
Other examples of classes of mappigns which can be described as solutions
of~\eqref{eq:osnovnoe-uravnenie}
with some 
functions~$F$
and~$G$
can be found in
\cite{BezrDK1987,Dair1986,Dair1992b,Dair1993b,Dair1995,Kopy1983,Kopy1984,Kopy1990,Kopy1995,Soko1991a,Soko1991b}.
The author has obtained some results on other properties of mappings of classes 
$\mathfrak G(K)$
and 
$\mathfrak F$
in
\cite{Egor2003,Egor2005,Egor2007,Egor2008,Egor2012,Egor2013,Egor2014a,Egor2014b}.

\section{Notation and Terminology}
Let 
$A$                                                                                      
be a set 
in~$\Bbb R^n$. 
The topological boundary 
of~$A$ 
is denoted 
by~$\partial A$. 
The diameter 
of~$A$ 
is defined as 
$\operatorname{diam}A:=\sup\{|x-y|:x,y\in A\}$. 
The outer Lebesgue measure 
of~$A$ 
is denoted 
by~$|A|$.

The set
$\Bbb R^{m\times n}:=\{\zeta=(\zeta_{\mu\nu})
_{\genfrac{}{}{0pt}{}{\mu=1,\dots,m}{\nu=1,\dots,n}}
:\zeta_{\mu\nu}\in\Bbb R,\ \mu=1,\dots,m,\ \nu=1,\dots,n\}$
consists of all real
($m\times n$)-matrices. 
We identify a matrix
$\zeta=(\zeta_{\mu\nu})
_{\genfrac{}{}{0pt}{}{\mu=1,\dots,m}{\nu=1,\dots,n}}
\in\Bbb R^{m\times n}$
with the linear mapping
$(\zeta_1,\dots,\zeta_m)\colon\Bbb R^n\to\Bbb R^m$, 
where
$\zeta_\mu(x):=\sum_{\nu=1}^n\zeta_{\mu\nu}x_\nu$,
$\mu=1,\dots,m$,
$x=(x_1,\dots,x_n)\in\Bbb R^n$.
The operator norm 
in~$\Bbb R^{m\times n}$
is defined as
$|\zeta|:=\sup\{|\zeta(x)|:x\in\Bbb R^n,\ |x|<1\}$.
The number of 
$k$-tuples 
of ordered indices in
$\Gamma_n^k:=\{I=(i_1,\dots,i_k):1\le i_1<\dots<i_k\le n,\
i_\varkappa\in\{1,\dots,n\},\ \varkappa=1,\dots,k\}$
equals the binomial coefficient
$\binom nk:=\frac{n!}{k!(n-k)!}$.
Given
$x\in\Bbb R^n$
and
$I\in\Gamma_n^k$,
we put
$x_I:=(x_{i_1},\dots,x_{i_k})\in\Bbb R^k$.
If
$I\in\Gamma_n^k$
and
$J\in\Gamma_m^k$,
then
$\det\nolimits_{JI}\zeta:=\det
\left(
\begin{smallmatrix}
\zeta_{j_1i_1}&\hdots&\zeta_{j_1i_k}\\
\vdots&\ddots&\vdots\\
\zeta_{j_ki_1}&\hdots&\zeta_{j_ki_k}
\end{smallmatrix}
\right)$
is the
$k\times k$-minors
of the matrix
$\zeta\in\Bbb R^{m\times n}$.

The Jacobi matrix of
$u=(u_1,\dots,u_m)\colon U\subset\Bbb R^n\to\Bbb R^m$
at a point
$x\in U$
is the matrix
$u'(x):=\bigl(\frac{\partial u_\mu}{\partial x_\nu}(x)\bigr)
_{\genfrac{}{}{0pt}{}{\mu=1,\dots,m}{\nu=1,\dots,n}}
$.
If
$I\in\Gamma_n^k$
and
$J\in\Gamma_m^k$
then 
$\frac{\partial u_J}{\partial x_I}(x)
=\frac{\partial(u_{j_1},\dots,u_{j_k})}
{\partial(x_{i_1},\dots,x_{i_k})}(x)
:=\det\nolimits_{JI}u'(x)$
and
$\partial_Iu_\mu(x)
:=\left(
\frac{\partial u_\mu}{\partial x_{i_1}}(x),\dots,\frac{\partial u_\mu}{\partial x_{i_1}}(x)
\right)$,
$\mu=1,\dots,m$.

Let
$\mathcal V$
be a real vector space
equipped with a 
norm~$|\cdot|$. 
We say that a function
$\Phi\colon\mathcal V\to\Bbb R$
is
{\it positively homogeneous of degree}
$p\in\Bbb R$
if
$\Phi(tx)=t^p\Phi(x)$
for all
$t>0$
and
$x\in\mathcal V\setminus\{0\}$.

Following 
Ch.~B.~Morrey~\cite{Morr1966},
we say that a continuous function
$F\colon\Bbb R^{m\times n}\to\Bbb R$
is 
{\it quasiconvex}, 
if
\begin{equation}
\label{eq:qc}
|B(0,1)|F(\zeta)\le\int_{B(0,1)}F(\zeta+\varphi'(x))\,dx
\end{equation}
for all
$\varphi\in C_0^\infty(B(0,1);\Bbb R^m)$
and
$\zeta\in\Bbb R^{m\times n}$.
Let 
$p\ge1$. 
Following 
M.~A.~Sychev~\cite{Sych1998}, 
we say that a quasiconvex 
function~$F$ 
is {\it strictly 
$p$-quasiconvex} 
if, for 
$\zeta\in\Bbb R^{m\times n}$ 
and 
$\varepsilon,C>0$, 
there is 
$\delta=\delta(\zeta,\varepsilon,C)>0$ 
such that, for each mapping 
$\varphi\in C_0^\infty(B(0,1);\Bbb R^m)$ 
satisfying 
$\|\varphi'\|_{L^p(B(0,1);\Bbb R^{m\times n})}
\le C|B(0,1)|^{1/p}$,
the condition
$\int_{B(0,1)}F(\zeta+\varphi'(x))\,dx
\le|B(0,1)|(F(\zeta)+\delta)$
implies 
$|\{x\in B(0,1):|\varphi'(x)|\ge\varepsilon\}|
\le\varepsilon|B(0,1)|$.
Observe that in the mathematical literature the term strictly 
quasiconvexity is also used for another property (which is close 
but nonequivalent to ours) consisting in the fact that the strict 
inequality in the definition of 
quasiconvexity~\eqref{eq:qc} 
is valid for nonzero 
mappings$~\varphi$ 
(for example, 
see~\cite{KnopS1984}). 
In this article we use the term in the sense of M.~A.~Sychev's 
definition~\cite{Sych1998}. 
In the case 
$p>1$ 
the notion of strictly 
$p$-quasiconvexity 
for 
functions~$F$ 
of this article is equivalent to the notion of strictly closed 
$p$-quasiconvexity 
from J.~Kristensen's
article~\cite{Kris1994} 
which is defined in terms of the theory of gradient Young measures 
(see 
\cite[Proposition 3.4]{Kris1994}).
Observe that we can replace the ball
$B(0,1)$ 
in the definitions of quasiconvexity and strictly 
$p$-quasiconvexity
by an arbitrary bounded 
domain~$U$ 
with 
$|\partial U|=0$
(for example,
see~\cite{Mull1999}).
A function
$G\colon\Bbb R^{m\times n}\to\Bbb R$
is a 
{\it null Lagrangian}
if both 
functions~$G$ 
and~$-G$ 
are quasiconvex. 
The term ``null Lagrangian'' appeared due to the following fact: 
The Euler--Lagrange equation corresponding to the variational 
integral
$\int_UG(u'(x))\,dx$
with null 
Lagrangian~$G$
holds identically for all admissible mappings
$u\colon U\subset\Bbb R^n\to\Bbb R^m$
(see~\cite{Ball1977b} 
and also~\cite{BallCO1981,IwanM2001,Mull1999}).
The only the affine combinations of minors (called 
{\it quasiaffine functions}) 
are null Lagrangians 
\cite{Edel1969,Land1942} 
(also see
\cite{Ball1977a,Ball1977b,BallCO1981,IwanM2001,Morr1966,Mull1999}); 
i.e.
\begin{equation}
\label{eq:reprsntn-nl}
G(\zeta)=\gamma_0
+\sum_{k=1}^{\min\{m,n\}}\sum_{J\in\Gamma_m^k,I\in\Gamma_n^k}
\gamma_{JI}\det\nolimits_{JI}\zeta,\quad 
\zeta\in\Bbb R^{m\times n},
\end{equation}
for some
$\gamma_0,\gamma_{JI}\in\Bbb R$.

Let 
$C_{\operatorname{loc}}(V;\Bbb R^m)$
be the space 
$C(V;\Bbb R^m)$
furnished with the topology of locally uniform convergence.

\section{Statement of the Main Results}
\label{sc:main-reslts}

Let 
$n,m,k\in\Bbb N$
such that
$2\le k\le\min\{n,m\}$.
We need the following hypothesis on continuous functions
$F\colon\Bbb R^{m\times n}\to\Bbb R$
and 
$G\colon\Bbb R^{m\times n}\to\Bbb R$
(see~\cite{Egor2008}):

(H1)
$F$
is quasiconvex;

(H1$'$) 
$F$ 
is strictly 
$k$-quasiconvex;

(H2)
$G$
is a null Lagrangian;
 
(H3)
$F$
and~$G$
are positively homogeneous of 
degree~$k$; 

(H4)
$\sup\{K\ge0:F(\zeta)\ge KG(\zeta),\ 
\zeta\in\Bbb R^{m\times n}\}=1$;

(H5)
$c_F:=\inf\{F(\zeta):\zeta\in\Bbb R^{m\times n},\ 
|\zeta|=1\}>0$.

By~(H3), the 
representation~\eqref{eq:reprsntn-nl} 
for the null 
Lagrangian~$G$
consists only of 
($k\times k$)-minors; 
i.e.,
\begin{equation}
\label{eq:reprsntn-G}
G(\zeta)=\sum_{J\in\Gamma_m^k,I\in\Gamma_n^k}
\gamma_{JI}\operatorname{det}_{JI}\zeta,
\quad\zeta\in\Bbb R^{m\times n}.
\end{equation}
It follows 
from~(H4)
that 
$\mathfrak G=\mathfrak G(1)$. 

The following theorems are the main results of the present paper.

\begin{theorem}
\label{th:C-sredn-ustoichivost-podoblast} 
Suppose 
$F$ 
and~$G$ 
satisfy  
{\rm(H1)--(H5)}.
Let 
$V$ 
be a bounded domain 
in~$\Bbb R^n$,
$K\ge1$,
and let 
$\mathcal S\subset\mathfrak G(K)\cap C(V;\Bbb R^m)$
such that 
$\mathcal S$
is precompact in 
$C_{\operatorname{loc}}(V;\Bbb R^m)$.
Then for a compact subset 
$U\subset V$
there exists a function 
$\alpha(\varepsilon)=\alpha_{\mathcal S,U}(\varepsilon)$, 
$0\le\varepsilon<\varepsilon_0$, 
$\lim_{\varepsilon\to0}\alpha(\varepsilon)=\alpha(0)=0$, 
such that, for every mapping
$v\in\mathcal S$
with
$\|K(\cdot,v)-1\|_{L^1(V)}<\varepsilon_0$,
there is a mapping
$u\colon V\to\mathbb R^m$ 
in the 
class~$\mathfrak G$
such that
\begin{equation}
\label{eq:C-otsenka}
\|v-u\|_{C(U;\Bbb R^m)}\le\alpha(\|K(\cdot,v)-1\|_{L^1(V)}).
\end{equation}
\end{theorem}

\begin{theorem}
\label{th:W-sredn-ustoichivost-podoblast} 
Suppose 
$F$ 
and~$G$ 
satisfy {\rm(H1$'$)}
and {\rm(H2)--(H5)}.
Let 
$V$ 
be a bounded domain 
in~$\Bbb R^n$,
$K\ge1$,
and let 
$\mathcal S\subset\mathfrak G(K)\cap C(V;\Bbb R^m)$
such that 
$\mathcal S$
is precompact in 
$C_{\operatorname{loc}}(V;\Bbb R^m)$.
Then for a compact subset 
$U\subset V$
there exists a function 
$\beta(\varepsilon)=\beta_{\mathcal S,U}(\varepsilon)$, 
$0\le\varepsilon<\varepsilon_0$, 
$\lim_{\varepsilon\to0}\beta(\varepsilon)=\beta(0)=0$, 
such that, for every mapping
$v\in\mathcal S$
with
$\|K(\cdot,v)-1\|_{L^1(V)}<\varepsilon_0$,
there is a mapping
$u\colon V\to\mathbb R^m$ 
in the 
class~$\mathfrak G$
such that
\begin{equation}
\label{eq:W-otsenka}
\|v-u\|_{C(U;\Bbb R^m)}+\|v'-u'\|_{L^k(U;\Bbb R^{m\times n})}
\le\beta(\|K(\cdot,v)-1\|_{L^1(V)}).
\end{equation}
\end{theorem}

\section{Proof of Theorem~\ref{th:C-sredn-ustoichivost-podoblast}} 

To prove 
Theorem~\ref{th:C-sredn-ustoichivost-podoblast}
we need the following auxillary lemma 
from~\cite{Egor2008}.

\begin{lemma}[\hbox{\cite[Lemma~1]{Egor2008}}]
\label{l:ravnomernaya-ogranichennost}
Let
$F$
and~$G$
satisfy
{\rm(H2)--(H5)}. 
Let 
$K\ge1$,
$V\subset\mathbb R^n$
be a domain, and 
$\mathcal S=\{v\colon V\to\mathbb R^m\}\subset\mathfrak G(K)$.
Suppose that 
$\mathcal S$
is uniformly bounded in
$L_{\operatorname{loc}}^k(V;\mathbb R^m)$.
Then
$\mathcal S$
is uniformly bounded in
$W_{\operatorname{loc}}^{1,k}(V;\mathbb R^m)$.
\end{lemma}

Let us prove 
Theorem~\ref{th:C-sredn-ustoichivost-podoblast}.
Proceeding by way of contradiction, assume that there are a compact subset 
$U\subset V$,
a number
$\varepsilon>0$,
and a sequence 
$(v_l\in\mathcal S)$
with
$\|K(\cdot,v_l)-1\|_{L^1(V)}\le1/l$
such that the inequality
\begin{equation}
\label{eq:C-neravenstvo-posledovatelnosti}
\|v_l-u\|_{C(U;\Bbb R^m)}>\varepsilon 
\end{equation}
holds for all mappings 
$u:V\to\Bbb R^m$
of the 
class~$\mathfrak G$.
Since
$\mathcal S$
is precompact in 
$C_{\operatorname{loc}}(V;\Bbb R^m)$
and
$\|K(\cdot,v_l)-1\|_{L^1(V)}\to 0$
as
$l\to\infty$,
from the sequence 
$(v_l)$
we can extract subsequence 
(we denote it by 
$(v_l)$
again)
such that it converges locally uniformly 
in~$V$
to some mapping
$v\colon V\to\Bbb R^m$
and
\begin{equation}
\label{eq:predel-K-posledovatelnost}
K(\cdot,v_l)\to1\quad\text{a.e. in }V
\end{equation}
as
$l\to\infty$.
Since
$\mathcal S\subset\mathfrak G(K)$,
from 
Lemma~\ref{l:ravnomernaya-ogranichennost}
we obtain that the 
sequence~$(v_l)$
is uniformly bounded in
$W_{\operatorname{loc}}^{1,k}(V;\mathbb R^m)$.
It follows from the general properties of the Sobolev spaces that
$v\in W_{\operatorname{loc}}^{1,k}(V;\mathbb R^m)$
(for example, see
\cite[Chapter~I, Theorem~1.1]{Resh1982a}).
We have
\begin{equation}
\label{eq:nl-neravenstvo-posledovatelnost}
F(v_l'(x))\le K(x,v_l)G(v_l'(x))\quad\text{a.e. in }V.
\end{equation}
and
\begin{equation}
\label{eq:ravnomernaya-otsenka-K-posledovatelnost}
K(x,v_l)\le K\quad\text{a.e. in }V.
\end{equation}
Multiply both sides 
of~\eqref{eq:nl-neravenstvo-posledovatelnost}  
by an arbitrary nonnegative function 
$\eta\in C_0^\infty(V)$
and integrate 
over~$V$. 
Eventually, we obtain
$\int_V\eta F(v_l')\le\int_V\eta K(\cdot,v_l)G(v_l')$.
Passing to the limit in the last inequality 
over~$l$
and using 
the theorem on weak semicontinuity of the functionals of calculus of variations 
\cite[Theorem~II.4]{AcerF1984},
the theorem on weak continuity of minors 
\cite[Chapter~II, Lemma~4.9]{Resh1982a}, 
\eqref{eq:predel-K-posledovatelnost},
and~\eqref{eq:ravnomernaya-otsenka-K-posledovatelnost},
we obtain
\begin{multline}
\label{eq:neravenstvo-polunepreryvnost}
\int_V\eta F(v')\le\liminf_{l\to\infty}\int_V\eta F(v_l')\\
\le\limsup_{l\to\infty}\int_V\eta F(v_l')
\le\limsup_{l\to\infty}\int_V\eta K(\cdot,v_l)G(v_l')\le\int_V\eta G(v').
\end{multline}
By the arbitrariness 
of~$\eta$, 
the last inequality means validity 
of~\eqref{eq:nl-neravenstvo} 
for~$v$ 
with
$K=1$.
It follows that
$v\in\mathfrak G$.
The sequence 
$(v_l)$
converges locally uniformly 
in~$V$
to~$v$.
This contradicts the 
assumption~\eqref{eq:C-neravenstvo-posledovatelnosti}.
Theorem~\ref{th:C-sredn-ustoichivost-podoblast} 
is proven.

\section{Proof of Theorem~\ref{th:W-sredn-ustoichivost-podoblast}} 

To prove 
Theorem~\ref{th:W-sredn-ustoichivost-podoblast}
we need the following auxillary proposition 
from~\cite{Egor2008}.

\begin{proposition}[\hbox{\cite[Proposition~1]{Egor2008}}]
Let 
$p>1$, 
and suppose that 
$F:\Bbb R^{m\times n}\to\Bbb R$ 
is a strictly 
$p$-quasiconvex function satisfying
$c|\zeta|^p\le F(\zeta)\le C(|\zeta|^p+1)$,
$\zeta\in\Bbb R^{m\times n}$,
with some constants 
$0<c<C<\infty$. 
Let 
$V\subset\Bbb R^n$
be a bounded domain with Lipschitz boundary, and let
$(v_l)_{l\in\Bbb N}$, 
$v_l\in W^{1,p}(V;\Bbb R^m)$, 
be a sequence of mappings such that 
$v_l\to v$ 
in 
$L^1(V;\Bbb R^m)$ 
for some mapping
$v\in W^{1,p}(V;\Bbb R^m)$. 
Suppose that 
$\int_V\eta F(v_l')\to\int_V\eta F(v')<\infty$.
Then 
$v_l\to v$ 
in 
$W^{1,p}(V;\Bbb R^m)$.
\end{proposition}

Let us prove 
Theorem~\ref{th:W-sredn-ustoichivost-podoblast}.
Assume that thare is no function with nessary properties.
Then for some number
$\varepsilon>0$
and some compact subset 
$U\subset V$
and every 
$l\in\Bbb N$
there exists a mapping 
$v_l\in \mathcal S$
with
$\|K(\cdot,v_l)-1\|_{L^1(V)}\le1/l$
such that the inequality
\begin{equation}
\label{eq:W-neravenstvo-posledovatelnosti}
\|v_l-u\|_{C(U;\Bbb R^m)}+\|v_l'-u'\|_{L^k(U;\Bbb R^{m\times n})}>\varepsilon
\end{equation}
holds for each mapping 
$u\in\mathfrak G$.
Arguing as the proof of 
Theorem~\ref{th:C-sredn-ustoichivost-podoblast},
we obtain that for the sequence 
$(v_l)$
there is a subsequence 
(denote it again by 
$(v_l)$)
converging locally uniformly 
in~$V$
to some mapping
$v\colon V\to\Bbb R^m$
from the 
class~$\mathfrak G$
and 
satisfying~\eqref{eq:neravenstvo-polunepreryvnost}
for any nonnegative function 
$\eta\in C_0^\infty(V)$.
We have that
$v$
satisfies~\eqref{eq:osnovnoe-uravnenie}.
Combining~\eqref{eq:neravenstvo-polunepreryvnost}
with~\eqref{eq:osnovnoe-uravnenie},
we have
\begin{equation}
\label{eq:ravenstvo-nizhnego-verkhnego-predela}
\liminf_{l\to\infty}\int_V\eta F(v_l')
=\limsup_{l\to\infty}\int_V\eta F(v_l')
=\int_V\eta F(v').
\end{equation}
It means that there is a subsequence (denoted again 
by $(v_l)$) 
such that
$\int_V\eta F(v_l')\to\int_V\eta F(v')$.
Observe that this subsequence depends on the chosen 
function~$\eta$. 
Taking an appropriate collection 
of~$\eta$
and using 
Proposition~???, 
we find that there is a subsequence (for which we preserve the notation 
$(v_l)$)
such that
$\|v_l'-v'\|_{L^k(U;\Bbb R^{m\times n})}\to0$.
Using the locally uniform convergence of 
$(v_l)$ 
to 
$v\in \mathfrak G$, 
we arrive at a contradiction 
with~\eqref{eq:W-neravenstvo-posledovatelnosti} 
Theorem~\ref{th:W-sredn-ustoichivost-podoblast} 
is proven.

\bibliographystyle{amsunsrt}

\end{document}